\documentclass[leqno,11pt]{article}
\usepackage[spanish]{babel}
\usepackage[utf8]{inputenc}
\usepackage{amscd}
\usepackage{amsthm}
\usepackage{amssymb}
\usepackage{amsmath}
\usepackage{graphics}
\usepackage{graphicx}
\usepackage{verbatim}
\usepackage{color}
\newtheorem{mydef}{Definition}
\newtheorem{myexa}{Example}

\newtheorem{mytheo}{Theorem}
\newtheorem{mylemma}{Lemma}

\newtheorem{myprop}{Proposition}

\newtheorem*{myprem}{{\em Remarks}}
\newtheorem*{myproo}{{\em Proof}}

\newcommand{\pt}{\mbox{$\succ$\hspace{-1ex}$\longrightarrow$}}

\decimalpoint

\begin{document}

\begin{center}
	{\sc The Bayes Estimator of a Conditional Density: Consistency}\vspace{2ex}\\
	A.G. Nogales\vspace{2ex}\\
	Dpto. de Matem\'aticas, Universidad de Extremadura\\
	Avda. de Elvas, s/n, 06006--Badajoz, SPAIN.\\
	e-mail: nogales@unex.es
\end{center}
\vspace{.4cm}
\begin{quote}
	\hspace{\parindent} {\small {\sc Abstract.}  In a Bayesian framework we prove that the optimal estimator of a conditional density is consistent. 
	}
\end{quote}

\vspace{3ex}
\begin{itemize}
	\item[] \hspace*{-1cm} {\em AMS Subject Class.} (2020): 62F15, 62G07
	\item[] \hspace*{-1cm} {\em Key words and phrases: } Bayesian estimation of a conditional density, posterior predictive distribution.
\end{itemize}

\vspace{3ex}

\section{Introduction.}

In the context of Functional Data Analysis, the problems of estimating a density, conditional or not, or a regression curve, both from the frequentist point of view and from the Bayesian point of view, are of particular interest. Nogales (2022a) addresses the density estimation problem from a Bayesian point of view and obtains the Bayes estimator for certain loss functions that are derived naturally from the quadratic error loss function when trying to estimate real functions of the parameter; under mild regularity conditions, that Bayesian estimator is concretely proven to be the posterior predictive density, and in Nogales (2022c) this optimal estimator is proven to be consistent. Nogales (2022b) obtains, among other results, the  Bayes estimator of the conditional density, and it is  the main aim of this paper to show its consistency.

The works mentioned above include references to the posterior predictive distribution, its usefulness in Predictive Inference and other areas of Statistical Inference, and its computation. Among those references we want to highlight here Geisser (1993) and Gelman et al. (2014). Ghosal et al. (2017) includes powerful tools to achieve this end. 

Some examples are included to illustrate the main result of the paper. For ease of reading we reproduce here an appendix from Nogales (2022a) to recall  the basic concepts of Bayesian inference but, mainly, to explain the (somehow unusual) notation and terminology used in this work.

\section{The framework.}

Let $(\Omega,\mathcal A,\{P_\theta\colon\theta\in(\Theta,\mathcal T,Q)\})$ be a Bayesian statistical experiment and $X_i:(\Omega,\mathcal A,\{P_\theta\colon\theta\in(\Theta,\mathcal T,Q)\})\rightarrow(\Omega_i,\mathcal A_i)$, $i=1,2$, two statistics. Consider the Bayesian experiment image of $(X_1,X_2)$
$$(\Omega_1\times\Omega_2,\mathcal A_1\times\mathcal  A_2,\{P_\theta^{(X_1,X_2)}\colon\theta\in(\Theta,\mathcal T,Q)\}).
$$
In the next, we shall suppose that $P^{(X_1,X_2)}(\theta,A_{12}):=P_\theta^{(X_1,X_2)}(A_{12})$, $\theta\in\Theta, \  A_{12}\in \mathcal A_1\times\mathcal  A_2$, is a Markov kernel. 
Let us write $R_\theta=P_\theta^{(X_1,X_2)}$ and $p_j(x):=x_j$ for $j=1,2$, $x:=(x_1,x_2)\in\Omega_1\times\Omega_2$. Hence
 \begin{gather*}
 	P_\theta^{X_1}=R_\theta^{p_1},\quad P_\theta^{X_2|X_1=x_1}=R_\theta^{p_2|p_1=x_1}\quad\hbox{and}\quad
 	E_{P_\theta}(X_2|X_1=x_1)=E_{R_\theta}(p_2|p_1=x_1).
 \end{gather*}

Given an integer $n$, for $m=n$ (resp. $m=\mathbb N$),  the Bayesian experiment corresponding to a $n$-sized sample (resp. an infinite sample) of the joint distribution of $(X_1,X_2)$ is 
$$\big((\Omega_1\times\Omega_2)^{m},(\mathcal A_1\times\mathcal  A_2)^{m},\big\{R_\theta^{m}\colon\theta\in(\Theta,\mathcal T,Q)\big\}\big)\qquad \hbox{(1)}
$$
We write $R^{m}(\theta,A'_{12,m})=R_\theta^{m}(A'_{12,m})$ for $A'_{12,m}\in (\mathcal A_1\times\mathcal  A_2)^{m}$ and 
$$\Pi_{12,{m}}:=Q\otimes R^{m},
$$
for the joint distribution of the parameter and the sample, i.e.
$$\Pi_{12,{m}}(A'_{12,m}\times T)=\int_TR_\theta^{m}(A'_{12,m})dQ(\theta),\quad A'_{12,m}\in(\mathcal A_1\times\mathcal  A_2)^{m}, T\in\mathcal T.
$$
The corresponding prior predictive distribution $\beta_{12,{m},Q}^*$ on $(\Omega_1\times\Omega_2)^{m}$ is
$$\beta_{12,{m},Q}^*(A'_{12,m})=\int_\Theta R_\theta^{m}(A'_{12,m})dQ(\theta),\quad A'_{12,m}\in(\mathcal A_1\times\mathcal  A_2)^{m}.
$$
The posterior distribution is a Markov kernel $$R_{m}^*:((\Omega_1\times\Omega_2)^{m},(\mathcal A_1\times\mathcal A_2)^{m})\pt (\Theta,\mathcal T)
$$
such that, for all $A'_{12,m}\in(\mathcal A_1\times\mathcal  A_2)^{m}$ and $T\in\mathcal T$,
$$\Pi_{12,{m}}(A'_{12,m}\times T)=\int_TR_\theta^{m}(A'_{12,m})dQ(\theta)
=\int_{A'_{12,m}}R_{m}^*(x',T)d\beta_{12,{m},Q}^*(x').
$$
Let us write $R_{{m},x'}^*(T):=R_{m}^*(x',T)$.

The posterior predictive distribution on $\mathcal A_1\times\mathcal A_2$  is the Markov kernel
$${R_{{m}}^*}^{R}:((\Omega_1\times\Omega_2)^{m},(\mathcal A_1\times\mathcal  A_2)^m)\pt(\Omega_1\times\Omega_2,\mathcal A_1\times\mathcal  A_2)
$$ 
defined, for $x'\in(\Omega_1\times\Omega_2)^{m}$, by
$${R_{{m}}^*}^{R}(x',A_{12}):=
\int_\Theta R_\theta(A_{12})dR_{{m},x'}^*(\theta)
$$
It follows that, with obvious notations, 
$$\int_{\Omega_1\times\Omega_2}f(x)d{R_{{m},x'}^*}^{\!\!\!\!R}(x)=\int_\Theta\int_{\Omega_1\times\Omega_2}f(x)dR_\theta(x)dR_{{m},x'}^*(\theta)
$$
for any real random variable (r.r.v. for short) $f$ whose integral exists. 

We can also consider the posterior predictive distribution on $(\mathcal A_1\times\mathcal A_2)^{m}$  defined as the Markov kernel
$${R_{{m}}^*}^{R^{m}}:((\Omega_1\times\Omega_2)^{m},(\mathcal A_1\times\mathcal  A_2)^{m})\pt((\Omega_1\times\Omega_2)^{m},(\mathcal A_1\times\mathcal  A_2)^{m})
$$ 
such that 
$${R_{{m}}^*}^{R^{m}}(x',A'_{12,m}):=
\int_\Theta R_\theta^{m}(A'_{12,m})dR_{{m},x'}^*(\theta)
$$

We introduce some notations for $(x',x,\theta)\in(\Omega_1\times\Omega_2)^{m}\times(\Omega_1\times\Omega_2)\times\Theta$:
\begin{gather*}
\pi'_{m}(x',x,\theta):=x',\quad \pi_{m}(x',x,\theta):=x,\quad \pi_{j,m}(x',x,\theta):=x_j,\;\; j=1,2,\quad
q_{m}(x',x,\theta):=\theta\\ \pi'_{i,m}(x',x,\theta):=x'_i:=(x'_{i1},x'_{i2}),\quad \pi'_{(i),m}(x',x,\theta):=(x'_1,\dots,x'_i),
\end{gather*}
for $1\le i\le m$ (read $i\in\mathbb N$  if  $m=\mathbb N$). 

Let us consider the probability space  
\begin{gather*}
((\Omega_1\times\Omega_2)^m\times(\Omega_1\times\Omega_2)\times\Theta,(\mathcal A_1\times\mathcal A_2)^m\times(\mathcal A_1\times\mathcal A_2)\times\mathcal T,\Pi_m),\qquad \hbox{(2)}
\end{gather*}
where
$$\Pi_m(A'_{12,m}\times A_{12}\times T)=\int_T R_\theta(A_{12}) R_\theta^m(A'_{12,m})dQ(\theta),
$$
when $A'_{12,m}\in(\mathcal A_1\times\mathcal A_2)^m$, $A_{12}\in \mathcal A_1\times\mathcal A_2$ and $T\in\mathcal T$. 

So, for a r.r.v. $f$ on $((\Omega_1\times\Omega_2)^m\times(\Omega_1\times\Omega_2)\times\Theta,(\mathcal A_1\times\mathcal A_2)^m\times(\mathcal A_1\times\mathcal A_2)\times\mathcal T)$,
$$\int f d\Pi_m=\int_\Theta\int_{(\Omega_1\times\Omega_2)^{m}}
\int_{\Omega_1\times\Omega_2} f(x',x,\theta)dR_\theta(x)dR_\theta^m(x')dQ(\theta)\qquad {\rm (3)}
$$
provided that the integral exists. 
Moreover, for a r.r.v. $h$ on $((\Omega_1\times\Omega_2)\times\Theta,(\mathcal A_1\times\mathcal A_2)\times\mathcal T)$,
$$\int h d\Pi_m=\int_\Theta\int_{\Omega_1\times\Omega_2}h(x,\theta)dR_\theta(x)dQ(\theta)
=\int_{\Omega_1\times\Omega_2}\int_\Theta h(x,\theta)dR^*_{1,x}(\theta)d\beta^*_{12,1,Q}(x).
$$ 

The following proposition is straightforward.

\begin{myprop}\label{prop1}\rm 
	For $n\in\mathbb N$, we have that
	\begin{gather*}
\Pi_{\mathbb N}^{(\pi'_{(n),\mathbb N},\pi_{1,\mathbb N},q_{\mathbb N})}=\Pi_n,\qquad
\Pi_{\mathbb N}^{(\pi'_{(n),\mathbb N},\pi_{1,\mathbb N})}=\Pi_n^{(\pi'_{(n),n},\pi_{1,n})},\\
\Pi_m^{q_m}=Q,\quad \Pi_m^{(\pi'_m,q_m)}=\Pi_{12,m},\quad \Pi_m^{\pi'_m}=\beta_{12,m,Q}^*,\quad 
\Pi_m^{(\pi_m,q_m)}=\Pi_{12,1},\quad
\Pi_m^{\pi_m}=\beta_{12,1,Q}^*\\ \Pi_m^{\pi'_m|q_m=\theta}=R_\theta^m,\quad \Pi_m^{\pi_m|q_m=\theta}=R_\theta,\quad \Pi_m^{q_m|\pi'_m=x'}=R^*_{m,x'},\quad \Pi_m^{q_m|\pi_m=x}=R^*_{1,x'}
\end{gather*}\end{myprop}

In particular, the probability space (2) contains all the basic ingredients of the Bayesian experiment (1), i.e., the prior distribution, the sampling probabilities, the posterior distributions and the prior predictive distribution. When $m=\mathbb N$ it becomes the natural framework to address the problem considered in this paper, as we can see in the next.

\section{Consistency of the Bayes estimator of the conditional density}

In the finite case ($m=n$), according to Theorem 1 of Nogales (2022b), when $\mathcal A_2$ is separable, 
the conditional distribution of $p_2$ given $p_1=x_1$ with respect to the posterior predictive distribution ${R^*_{n,x'}}^{\!\!\!\!R}$
is the Bayes estimator of the conditional distribution   $R^{p_2|p_1}$ for the squared total variation loss function, i.e.,
\begin{gather*}
	\int_{(\Omega_1\times\Omega_2)^{n+1}\times\Theta}
	\sup_{A_2\in\mathcal A_2}\left|\left({R^*_{n,x'}}^{\!\!\!\!R}\right)^{p_2|p_1=x_1}-R_\theta^{p_2|p_1=x_1}(A_2)\right|^2
	d\Pi_n(x',x,\theta)\le
	\\
	\int_{(\Omega_1\times\Omega_2)^{n+1}\times\Theta}
	\sup_{A_2\in\mathcal A_2}\big|M(x',x_1,A_2)-R_\theta^{p_2|p_1=x_1}(A_2)\big|^2
	d\Pi_n(x',x,\theta)
\end{gather*}
for any estimator $M$ of the conditional distribution $R^{p_2|p_1}$. Recall from Nogales (2022b) that an estimator of the conditional distribution $P_\theta^{X_2|X_1}$ from a $n$-sized sample of the joint distribution of $(X_1,X_2)$ is a Markov kernel 
$$M:((\Omega_1\times\Omega_2)^n\times\Omega_1,(\mathcal A_1\times\mathcal A_2)^n\times\mathcal A_1)\pt (\Omega_2,\mathcal A_2)
$$
so that, being observed $x'=((x'_{11},x'_{21}),\dots,(x'_{1n},x'_{2n}))\in (\Omega_1\times \Omega_2)^n$, $M(x',x_1,\cdot)$ is a probability measure on $\mathcal A_2$ that is considered as an estimation of the conditional distribution $P_\theta^{X_2|X_1=x_1}$ for a given $x_1\in\Omega_1$.

We wonder if the Bayes risk of this estimator goes to 0 when $n$ grows to $\infty$, i.e. if
$$\lim_n
	\int_{(\Omega_1\times\Omega_2)^{n+1}\times\Theta}
\sup_{A_2\in\mathcal A_2}\left|\left({R^*_{n,x'}}^{\!\!\!\!R}\right)^{p_2|p_1=x_1}(A_2)-R_\theta^{p_2|p_1=x_1}(A_2)\right|^2
d\Pi_n(x',x,\theta)=0.\qquad (4)
$$

When the joint distribution $R_\theta=P_\theta^{(X_1,X_2)}$ has a density $f_\theta$ respect to the product of two $\sigma$-finite measures $\mu_1$ and $\mu_2$ on $\mathcal A_1$ and $\mathcal A_2$, resp., the conditional density is
$$f_\theta^{X_2|X_1=x_1}(x_2):=\frac{f_\theta(x_1,x_2)}{f_{\theta,X_1}(x_1)}
$$
for almost every $x_1$, where $f_{\theta,X_1}(x_1)$ stands for the marginal density of $X_1$.

According to Theorem 2 of Nogales (2022b), when $\mathcal A_2$ is separable, the Bayes estimator of the conditional density $f_\theta^{X_2|X_1}$ (from an $n$-sized sample) for the $L^1$-squared loss function is the  $\mu_2$-density $${f^*_{n,x'}}^{\!\!\!\!X_2|X_1=x_1}(x_2):=\frac{f^*_{n,x'}(x_1,x_2)}{f^*_{n,x',1}(x_1)}=\frac{\int_\Theta f_\theta(x_1,x_2) r^*_{n,x'}(\theta)dQ(\theta)}{\int_{\Omega_2}\int_\Theta f_\theta(x_1,t) r^*_{n,x'}(\theta)dQ(\theta)d\mu_2(t)}
$$ 
of the conditional distribution
$\left[{R^{*}_{n,x'}}^{\!\!\!\!R}\right]^{p_2|p_1}$ of $p_2$ given $p_1$ with respect to the posterior predictive distribution ${R^{*}_{n,x'}}^{\!\!\!\!R}$,
i.e.
\begin{gather*}
	\int_{(\Omega_1\times\Omega_2)^{n+1}\times\Theta}
	\left(\int_{\Omega_2}\left|{f^*_{n,x'}}^{\!\!\!\!X_2|X_1=x_1}(t)-f_\theta^{X_2|X_1=x_1}(t)\right|d\mu_2(t)\right)^2d\Pi_n(x',x,\theta)\le\\
	\int_{(\Omega_1\times\Omega_2)^{n+1}\times\Theta}
	\left(\int_{\Omega_2}\left|d(x',x_1,t)-f_\theta^{X_2|X_1=x_1}(t)\right|d\mu_2(t)\right)^2d\Pi_n(x',x,\theta),
\end{gather*}
for any estimator $d$ of the conditional density.

Recall from Nogales (2022b) that an estimator of the conditional density $f_\theta^{X_2|X_1}$ from a $n$-sized sample of the joint distribution of $(X_1,X_2)$ is a map 
$$d:((\Omega_1\times\Omega_2)^n\times\Omega_1\times\Omega_2,(\mathcal A_1\times\mathcal A_2)^n\times\mathcal A_1\times\mathcal A_2)\rightarrow (\mathbb R,\mathcal R)
$$
so that, being observed $x'=((x'_{11},x'_{21}),\dots,(x'_{1n},x'_{2n}))\in (\Omega_1\times \Omega_2)^n$, $d(x',x_1,\cdot)$  is considered as an estimation of the conditional density $f_\theta^{X_2|X_1=x_1}$ of $X_2$ given $X_1=x_1$.

In terms of densities, question (4) reads as follows:
$$\lim_n	\int_{(\Omega_1\times\Omega_2)^{n+1}\times\Theta}
\left(\int_{\Omega_2}\left|{f^*_{n,x'}}^{\!\!\!\!X_2|X_1=x_1}(t)-f_\theta^{X_2|X_1=x_1}(t)\right|d\mu_2(t)\right)^2d\Pi_n(x',x,\theta)=0.\qquad (5)$$

Without loss of generality we can suppose that $\mu_2$ is a probability measure. 

In the next we will assume the following additional regularity conditions: 
\begin{itemize}
	\item[(i)] $(\Omega_i,\mathcal A_i)$, $i=1,2$, are  standard Borel spaces, 
	
	\item[(ii)] $\Theta$ is a Borel subset of a Polish space and  $\mathcal T$ is its Borel $\sigma$-field, and
	\item[(iii)] $\{R_\theta\colon \theta\in\Theta\}$ is identifiable.
\end{itemize}

Let us consider the auxiliary Bayesian experiment
$$(\Omega_2\times(\Omega_1\times\Omega_2)^{\mathbb N}\times(\Omega_1\times\Omega_2)\times\Theta,\mathcal A_2\times (\mathcal A_1\times\mathcal A_2)^{\mathbb N}\times(\mathcal A_1\times\mathcal A_2)\times\mathcal T,\mu_2\times \Pi_{\mathbb N}).
$$
For $t\in\Omega_2$, $x'\in(\Omega_1\times\Omega_2)^{\mathbb N}$, $x\in\Omega_1\times\Omega_2$, $\theta\in\Theta$ and $n\in\mathbb N$, we write
\begin{gather*}
I'(t,x',x,\theta)=x',\quad
I'_n(t,x',x,\theta):=x'_n:=(x'_{n1},x'_{n2}),\quad I'_{(n)}(t,x',x,\theta):=(x'_1,\dots,x'_n),\\
 I(t,x',x,\theta)=x,\quad I_j(t,x',x,\theta)=x_j,\; j=1,2,\quad K(t,x',x,\theta)=t,\quad J(t,x',x,\theta)=\theta. 
\end{gather*}

A similar result to Proposition 1 holds for these projections. 
The following lemma will be useful to solve our problem. 

\begin{mylemma}\label{lemma1}\rm Let $Y(t,x',x,\theta):=f_\theta^{X_2|X_1=x_1}(t)$.
	
(i)	For $n\in\mathbb N$  we have that
	\begin{gather*}\begin{split}
			{f^*_{n,x'_{(n)}}}^{\!\!\!\!X_2|X_1=x_1}(t)&=E_{\mu_2\times\Pi_{\mathbb N}}(Y|(K,I'_{(n)},I_1)=(t,x'_{(n)},x_1))\\
			&=\int_\Theta Y(t,x',x,\theta)d(\mu_2\times\Pi_{\mathbb N})^{J|(K,I'_{(n)},I_1)=(t,x'_{(n)},x_1)}(\theta).
	\end{split}\end{gather*}

(ii) $${f^*_{\mathbb N,x'}}^{\!\!\!\!X_2|X_1=x_1}(t)=E_{\mu_2\times\Pi_{\mathbb N}}(Y|(K,I',I_1)=(t,x',x_1)).
$$
\end{mylemma}

When $\mathcal A'_{(n)}:=(K,I'_{(n)},I_1)^{-1}(\mathcal A_2\times(\mathcal A_1\times\mathcal A_2)^n\times\mathcal A_1)$, 
we have that  $(\mathcal A'_{(n)})_n$ is an increasing sequence of  sub-$\sigma$-fields of $\mathcal A_2\times(\mathcal A_1\times\mathcal A_2)^{\mathbb N}\times \mathcal A_1$ such that  $\mathcal A_2\times(\mathcal A_1\times\mathcal A_2)^{\mathbb N}\times \mathcal A_1=\sigma(\cup_n\mathcal A'_{(n)})$. According to the martingale convergence theorem of Lévy, if $Y$ is $(\mathcal A_2\times(\mathcal A_1\times\mathcal A_2)^{\mathbb N}\times \mathcal A_1\times\mathcal T)$-measurable and $\mu_2\times\Pi_{\mathbb N}$-integrable, then 
$$E_{\mu_2\times \Pi_{\mathbb N}}(Y|\mathcal A'_{(n)})
$$
converges $(\mu_2\times\Pi_{\mathbb N})$-a.e. and in $L^1(\mu_2\times\Pi_{\mathbb N})$ to $E_{\mu_2\times\Pi_{\mathbb N}}(Y|\mathcal A_2\times(\mathcal A_1\times\mathcal A_2)^{\mathbb N}\times\mathcal A_1)$. 

Let us consider the $\mu_2\times\Pi_{\mathbb N}$-integrable function 
$$Y(t,x',x,\theta):=f_\theta^{X_2|X_1=x_1}(t). 
$$
Hence, it follows from the aforementioned theorem of Lévy that
$$\lim_n {f^*_{n,x'_{(n)}}}^{\!\!\!\!X_2|X_1=x_1}(t) =
{f^*_{\mathbb N,x'}}^{\!\!\!\!X_2|X_1=x_1}(t),\quad \mu_2\times\Pi_{\mathbb N}-\hbox{a.e.} \qquad (10)
$$
and
$$\lim_n\int_{\Omega_2\times(\Omega_1\times\Omega_2)^{\mathbb N}\times(\Omega_1\times\Omega_2)\times\Theta}\left|{f^*_{n,x'_{(n)}}}^{\!\!\!\!X_2|X_1=x_1}(t)-{f^*_{\mathbb N,x'}}^{\!\!\!\!X_2|X_1=x_1}(t)\right| d(\mu_2\times\Pi_{\mathbb N})(t,x',x,\theta)=0,
$$
i.e., 
$$\lim_n\int_{(\Omega_1\times\Omega_2)^{\mathbb N}\times(\Omega_1\times\Omega_2)\times\Theta}\int_{\Omega_2}\left|{f^*_{n,x'_{(n)}}}^{\!\!\!\!X_2|X_1=x_1}(t)-{f^*_{\mathbb N,x'}}^{\!\!\!\!X_2|X_1=x_1}(t)\right| d\mu_2(t)d\Pi_{\mathbb N}(x',x,\theta)=0.\qquad (11)
$$

On the other hand, as a consequence of a known theorem of Doob (see Theorem 6.9 and Proposition 6.10 from Ghosal et al. (2017, p. 129, 130)) we have that, for every $t\in\Omega_2$ and $x_1\in\Omega_1$, 
$$\lim_n \int_{\Theta}f_{\theta'}^{X_2|X_1=x_1}(t)
d(\mu_2\times\Pi_{\mathbb N})^{J|(K,I'_{(n)},I_1)=(t,x'_{(n)},x_1)}(\theta')=f_{\theta}^{X_2|X_1=x_1}(t)
,\quad R_\theta^{\mathbb N}-\hbox{a.e.}
$$
for $Q$-almost every $\theta$. Hence, according to Lemma \ref{lemma1} (i), 
$$\lim_n {f^*_{n,x'_{(n)}}}^{\!\!\!\!X_2|X_1=x_1}(t)=f_{\theta}^{X_2|X_1=x_1}(t)
,\quad R_\theta^{\mathbb N}-\hbox{a.e.}
$$
for $Q$-almost every $\theta$.
 In particular,
$$\lim_n {f^*_{n,x'_{(n)}}}^{\!\!\!\!X_2|X_1=x_1}(t)=f_{\theta}^{X_2|X_1=x_1}(t)
,\quad (\mu_2\times\Pi_{\mathbb N})-\hbox{a.e.}
$$
In this sense we can say that the posterior predictive conditional density ${f^*_{n,x'_{(n)}}}^{\!\!\!\!X_2|X_1=x_1}(t)$ is a strongly consistent estimator of the conditional density $f_{\theta}^{X_2|X_1=x_1}(t)$. 

From this and (10) we obtain that
$$ {f^*_{\mathbb N,x'}}^{\!\!\!\!X_2|X_1=x_1}(t)=f_{\theta}^{X_2|X_1=x_1}(t)
,\quad (\mu_2\times\Pi_{\mathbb N})-\hbox{a.e.}
$$
According to (11) we obtain that
$$\lim_n\int_{(\Omega_1\times\Omega_2)^{\mathbb N}\times(\Omega_1\times\Omega_2)\times\Theta}\int_{\Omega_2}\left|{f^*_{n,x'_{(n)}}}^{\!\!\!\!X_2|X_1=x_1}(t)-f_{\theta}^{X_2|X_1=x_1}(t)\right| d\mu_2(t)d\Pi_{\mathbb N}(x',x,\theta)=0,\qquad (12)
$$
which proves that the risk function of the Bayes estimator ${f^*_{n,x'_{(n)}}}^{\!\!\!\!X_2|X_1=x_1}(t)$ of the conditional density $f_{\theta}^{X_2|X_1=x_1}(t)$ converges to 0 for the $L^1$-loss function. 
Since
$$
\sup_{A_2\in\mathcal A_2}\left|\left({R^*_{n,x'}}^{\!\!\!\!R}\right)^{p_2|p_1=x_1}(A_2)-R_\theta^{p_2|p_1=x_1}(A_2)\right|=\frac12\int_{\Omega_2}\left|{f^*_{n,x'_{(n)}}}^{\!\!\!\!X_2|X_1=x_1}(t)-f_{\theta}^{X_2|X_1=x_1}(t)\right| d\mu_2(t)
$$
we have that
$$\lim_n\int_{(\Omega_1\times\Omega_2)^{\mathbb N}\times(\Omega_1\times\Omega_2)\times\Theta}\sup_{A_2\in\mathcal A_2}\left|\left({R^*_{n,x'}}^{\!\!\!\!R}\right)^{p_2|p_1=x_1}(A_2)-R_\theta^{p_2|p_1=x_1}(A_2)\right|d\Pi_{\mathbb N}(x',x,\theta)=0,\qquad (13)
$$
which proves that that the risk function of the Bayes estimator $\left({R^*_{n,x'}}^{\!\!\!\!R}\right)^{p_2|p_1=x_1}$ of the conditional distribution $R_\theta^{p_2|p_1=x_1}$ converges to 0 for the total variation loss function.

We wonder if these results remain true for the squared versions of the loss functions. The answer is yes because bounded $L^1$-convergence on a probability space implies $L^2$-convergence.

So, we have proved the following result. 

\begin{mytheo}\rm Let $(\Omega,\mathcal A,\{P_\theta\colon\theta\in(\Theta,\mathcal T,Q)\})$ be a Bayesian statistical experiment and $X_i:(\Omega,\mathcal A,\{P_\theta\colon\theta\in(\Theta,\mathcal T,Q)\})\rightarrow(\Omega_i,\mathcal A_i)$, $i=1,2$, two statistics. Consider the Bayesian experiment image of $(X_1,X_2)$
	$$(\Omega_1\times\Omega_2,\mathcal A_1\times\mathcal  A_2,\{R_\theta\colon\theta\in(\Theta,\mathcal T,Q)\}).
	$$
	where $R_\theta:=P_\theta^{(X_1,X_2)}$ is supposed dominated by the product $\mu_1\times\mu_2$ of two $sigma$-finite measures $\mu_1$ and $\mu_2$ on $\mathcal A_1$ and $\mathcal A_2$, resppectively. 
	Let us suppose that $(\Omega_1\times\Omega_2,\mathcal A_1\times\mathcal  A_2)$ is a  Borel standar space, that $\Theta$ is a Borel subset of a Polish space and $\mathcal T$ is its Borel  $\sigma$-field. Suppose also that the likelihood function  $\mathcal L(x_1,x_2;\theta):=f_\theta(x_1,x_2)=\frac{dR_\theta}{d(\mu_1\times\mu_2)}(x_1,x_2)$ is $\mathcal A_1\times\mathcal A_2\times\mathcal T$-measurable and the family $\{R_\theta\colon \theta\in\Theta\}$ is identifiable. Then:
	\begin{itemize}
		\item[(a)] The posterior predictive conditional density ${f^*_{n,x'_{(n)}}}^{\!\!\!\!X_2|X_1}$ 
		is the Bayes estimator of the conditional density  $f_\theta^{X_2|X_1}$ in the product experiment  $((\Omega_1\times\Omega_2)^n,(\mathcal A_1\times\mathcal  A_2)^n,\{R_\theta^n\colon \theta\in(\Theta,\mathcal T,Q)\})$ for the  squared $L^1$ loss function. Moreover the risk function converges to 0 both for the $L^1$ loss function and the squared  $L^1$ loss function. 
		
		\item[(b)] The posterior predictive conditional distribution  $\left({R^*_{n,x'}}^{\!\!\!\!R}\right)^{p_2|p_1=x_1}$ is the Bayes estimator of the sampling conditional distribution  $P_\theta^{X_2|X_1}$ in the product experiment  $((\Omega_1\times\Omega_2)^n,(\mathcal A_1\times\mathcal  A_2)^n,\{R_\theta^n\colon \theta\in(\Theta,\mathcal T,Q)\})$ for the squared variation total loss function. Moreover the risk function converges to 0 both for the total variation loss function and the squared total variation loss function. 
		
		\item[(c)] The posterior predictive density is a strongly consistent estimator of the conditional density $f_{\theta}^{X_2|X_1}$, in the sense that
		$$\lim_n {f^*_{n,x'_{(n)}}}^{\!\!\!\!X_2|X_1=x_1}(t)=f_{\theta}^{X_2|X_1=x_1}(t),\quad (\mu_2\times\Pi_{\mathbb N})-\hbox{a.e.}
		$$
		for $Q$-almost every $\theta\in\Theta$. 
		
	\end{itemize}
\end{mytheo}

\section{Proof of the Lemma}

\begin{myproo}\rm 
	(i)	According to Lemma 1 of Nogales (2022b), we have that, for all $A'_{12,n}\in(\mathcal A_1\times\mathcal A_2)^n$ and all $A_1\in\mathcal A_1$, 
	\begin{gather*}
		\int_{A'_{12,n}\times A_1\times\Omega_2\times \Theta}R_\theta^{p_2|p_1=x_1}(A_2)d\Pi_n(x',x,\theta)=\\
		\int_{A'_{12,n}\times A_1}\left[{R^{*}_{n,x'}}^{\!\!\!\!R}\right]^{p_2|p_1=x_1}\!\!(A_2)d{\Pi_n}^{(\pi',p_1)}(x',x_1).
	\end{gather*}
	From this and Proposition \ref{prop1} we obtain that, given $A_2\in\mathcal A_2$, 
	\begin{gather*}
		\int_{A_2\times(A'_{12,n}\times(\Omega_1\times\Omega_2)^{\mathbb N})\times(A_1\times\Omega)\times\Theta}f_\theta^{X_2|X_1=x_1}(t)
		d(\mu_2\times\Pi_{\mathbb N})(t,x',x,\theta)=\\
		\int_{A'_{12,n}\times A_1\times\Omega_2\times \Theta}R_\theta^{p_2|p_1=x_1}(A_2)d\Pi_n(x'_{(n)},x,\theta)=\\
		\int_{A'_{12,n}\times A_1}\left[{R^{*}_{n,x'_{(n)}}}^{\!\!\!\!R}\right]^{p_2|p_1=x_1}\!\!(A_2)d{\Pi_n}^{(\pi'_{(n),n},\pi_{1,n})}(x'_{(n)},x_1)=\\
		\int_{A_2\times A'_{12,n}\times A_1}
		{f^*_{n,x'_{(n)}}}^{\!\!\!\!X_2|X_1=x_1}(t)d(\mu_2\times\Pi_{\mathbb N})^{(K,I'_{(n)},I_1)}(t,x'_{(n)},x_1),
	\end{gather*}
	which finishes the proof of the first part of (i). Finally, desintegrating 
	$$\mu_2\times\Pi_{\mathbb N}=(\mu_2\times\Pi_{\mathbb N})^{(J|(K,I'_{(n)},I_1))}\otimes (\mu_2\times\Pi_{\mathbb N})^{(K,I'_{(n)},I_1)}
	$$
	we obtain 
	\begin{gather*}
		\int_{A_2\times(A'_{12,n}\times(\Omega_1\times\Omega_2)^{\mathbb N})\times(A_1\times\Omega)\times\Theta}f_\theta^{X_2|X_1=x_1}(t)
		d(\mu_2\times\Pi_{\mathbb N})(t,x',x,\theta)=\\
		\int_{A_2\times A'_{12,n}\times A_1}\int_{\Theta} f_\theta^{X_2|X_1=x_1}(t)
		d(\mu_2\times\Pi_{\mathbb N})^{J|(K,I'_{(n)},I_1)=(t,x'_{(n)},x_1)}(\theta)d(\mu_2\times\Pi_{\mathbb N})^{(K,I'_{(n)},I_1)}(t,x'_{(n)},x_1)	
	\end{gather*}
	which proves the second part of (i).

	(ii) An analogous reasoning would work in this case if we can prove that, for all $A'_{12,\mathbb N}\in(\mathcal A_1\times\mathcal A_2)^{\mathbb N}$ and all $A_1\in\mathcal A_1$, 
	\begin{gather*}
		\int_{A'_{12,\mathbb N}\times A_1\times\Omega_2\times \Theta}R_\theta^{p_2|p_1=x_1}(A_2)d\Pi_{\mathbb N}(x',x,\theta)=\\
		\int_{A'_{12,{\mathbb N}}\times A_1}\left[{R^{*}_{{\mathbb N},x'}}^{\!\!\!\!R}\right]^{p_2|p_1=x_1}\!\!(A_2)d{\Pi_{\mathbb N}}^{(\pi',p_1)}(x',x_1).\qquad\hbox{(6)}
		\	\end{gather*}
	According to (3), we have that, for all $A'_{12,{\mathbb N}}\in(\mathcal A_1\times\mathcal A_2)^{\mathbb N}$ and all $A_1\in\mathcal A_1$, 
	\begin{gather*}
		\int_{A'_{12,{\mathbb N}}\times A_1\times\Omega_2\times \Theta}R_\theta^{p_2|p_1=x_1}(A_2)d\Pi_{\mathbb N}(x',x,\theta)=\\
		\int_\Theta\int_{A'_{12,{\mathbb N}}}\int_{A_1}R_\theta^{p_2|p_1=x_1}(A_2)dR_\theta^{p_1}(x_1)dR_\theta^{\mathbb N}(x')dQ(\theta)=\\
		\int_\Theta R_\theta^{\mathbb N}(A'_{12,{\mathbb N}})R_\theta(A_1\times A_2)dQ(\theta).\quad \hbox{(7)}
	\end{gather*}

	Note that, by definition of conditional distribution,
	$${R^{*}_{{\mathbb N},x'}}^{\!\!\!\!R}(A_1\times A_2)=\int_{A_1}\left[{R^{*}_{{\mathbb N},x'}}^{\!\!\!\!R}\right]^{p_2|p_1=x_1}\!\!(A_2)d\!\left[{R^{*}_{{\mathbb N},x'}}^{\!\!\!\!R}\right]^{p_1}\!\!(x_1),\qquad \hbox{(8)}
	$$
	and, by definition of posterior predictive distribution,
	$${R^{*}_{{\mathbb N},x'}}^{\!\!\!\!R}(A_1\times A_2)=\int_{\Theta}R_\theta(A_1\times A_2)dR^*_{{\mathbb N},x'}(\theta).
	$$
	
	Notice that, for any $A'_{12,{\mathbb N}}\in(\mathcal A_1\times\mathcal A_2)^{\mathbb N}$, being $R^*_{{\mathbb N},x'}=\Pi_{\mathbb N}^{q_{\mathbb N}|\pi'_{\mathbb N}=x'}$, 
	\begin{gather*}
		\int_{A'_{12,{\mathbb N}}}\int_{\Theta} R_\theta(A_1\times A_2) dR^*_{{\mathbb N},x'}(\theta)d\Pi_{\mathbb N}^{\pi'_{\mathbb N}}(x')=
		\int_{A'_{12,{\mathbb N}}}\int_{\Theta}R_\theta(A_1\times A_2)d\Pi_{\mathbb N}^{q_{\mathbb N}|\pi'_{\mathbb N}=x'}(\theta)d\Pi_{\mathbb N}^{\pi'_{\mathbb N}}(x')=\\
		\int_{A'_{12,{\mathbb N}}\times\Theta}R_\theta(A_1\times A_2)d\Pi_{\mathbb N}^{(\pi'_{\mathbb N},q_{\mathbb N})}(x',\theta)=\int_{{\pi'_{\mathbb N}}^{-1}(A'_{12,{\mathbb N}})}R_\theta(A_1\times A_2)d\Pi_{\mathbb N}(x',x,\theta),
	\end{gather*}
	proving that
	$$\int_{\Theta}R_\theta(A_1\times A_2)dR^*_{{\mathbb N},x'}(\theta)=	
	E_{\Pi_{\mathbb N}}(r_{A_1\times A_2}|\pi'_{\mathbb N}=x')
	$$
	where $r_{A_1\times A_2}(\theta):=R_\theta(A_1\times A_2)$ and, hence
	$${R^{*}_{{\mathbb N},x'}}^{\!\!\!\!R}(A_1\times A_2)=E_{\Pi_{\mathbb N}}(r_{A_1\times A_2}|\pi'_{\mathbb N}=x').
	$$ 
	
	So, by definition of conditional expectation,
	$$\int_{A'_{12,{\mathbb N}}}{R^{*}_{{\mathbb N},x'}}^{\!\!\!\!R}(A_1\times A_2)d\Pi_{\mathbb N}^{\pi'_{\mathbb N}}(x')=\int_{A'_{12,{\mathbb N}}\!\times\Omega_1\times\Omega_2\times \Theta}R_\theta(A_1\times A_2)d\Pi_{\mathbb N}(x',x,\theta).\qquad \hbox{(9)}
	$$
	But the second term of this equation is
	\begin{gather*}
		\int_{A'_{12,{\mathbb N}}\times\Omega_1\times\Omega_2\times \Theta}R_\theta(A_1\times A_2)d\Pi_{\mathbb N}(x',x,\theta)=\\\int_\Theta\int_{A'_{12,{\mathbb N}}}R_\theta(A_1\times A_2)dR_\theta^{\mathbb N}(x')dQ(\theta)=\int_\Theta R_\theta^{\mathbb N}(A'_{12,{\mathbb N}})R_\theta(A_1\times A_2)dQ(\theta)
	\end{gather*}
	which coincides with (7). 
	
	From (6) and Proposition \ref{prop1} we obtain that, given $A_2\in\mathcal A_2$, 
	\begin{gather*}
		\int_{A_2\times(A'_{12,{\mathbb N}}\times(\Omega_1\times\Omega_2)^{\mathbb N}\times(A_1\times\Omega)\times\Theta}f_\theta^{X_2|X_1=x_1}(t)
		d(\mu_2\times\Pi_{\mathbb N})(t,x',x,\theta)=\\
		\int_{A'_{12,{\mathbb N}}\times A_1\times\Omega_2\times \Theta}R_\theta^{p_2|p_1=x_1}(A_2)d\Pi_{\mathbb N}(x',x,\theta)=\\
		\int_{A'_{12,{\mathbb N}}\times A_1}\left[{R^{*}_{{\mathbb N},x'}}^{\!\!\!\!R}\right]^{p_2|p_1=x_1}\!\!(A_2)d{\Pi_{\mathbb N}}^{(\pi'_{\mathbb N},\pi_{1,{\mathbb N}})}(x',x_1)=\\
		\int_{A_2\times A'_{12,{\mathbb N}}\times A_1}
		{f^*_{{\mathbb N},x'}}^{\!\!\!\!X_2|X_1=x_1}(t)d(\mu_2\times\Pi_{\mathbb N})^{(K,I',I_1)}(t,x',x_1),
	\end{gather*}
	which finishes the proof of (ii). $\Box$
\end{myproo}

\section{Examples.}\label{secex}

\begin{myexa}\rm Let us suppose that, for $\theta,\lambda,x_1>0$, $P_\theta^{X_1}= G(1,\theta^{-1})$, $P_\theta^{X_2|X_1=x_1}=G(1,(\theta x_1)^{-1})$ and $Q=G(1,\lambda^{-1})$, where $G(\alpha,\beta)$ denotes the gamma distribution of parameters $\alpha,\beta>0$. Hence the joint density of $X_1$ and $X_2$ is
	$$f_\theta(x_1,x_2)=\theta^2 x_1\exp\{-\theta x_1(1+x_2)\}I_{]0,\infty[^2}(x_1,x_2).
	$$
	It is shown in Nogales (2022b), Example 1, that 
the Bayes estimator of the conditional density $f_\theta^{X_2|X_1=x_1}(x_2)=\theta x_1\exp\{-\theta x_1x_2\}I_{]0,\infty[}(x_2)$ is, for $x_1,x_2>0$, 
\begin{gather*}
{f^*_{n,x'}}^{\!\!\!\!X_2|X_1=x_1}(x_2):=\frac{f^*_{n,x'}(x_1,x_2)}{f^*_{n,x',1}(x_1)}=\\(2n+2)x_1\frac{[\lambda+x_1+\sum_{i=1}^nx'_{i1}(1+x'_{i2})]^{2n+2}}{[\lambda+x_1(1+x_2)+\sum_{i=1}^nx'_{i1}(1+x'_{i2})]^{2n+3}}=\\
\frac{(2n+2)x_1a_n(x',x_1)^{2n+2}}{(x_1x_2+a_n(x',x_1))^{2n+3}}
\end{gather*}
where 
$$a_n(x',x_1)=\lambda+x_1+\sum_{i=1}^nx'_{i1}(1+x'_{i2}).
$$
Theorem 1 shows that this a strongly consistent estimator of de conditional density $f_\theta^{X_2|X_1=x_1}$.
	\end{myexa} 
	
\begin{myexa}\rm Let us suppose that $X_1$ has a Bernoulli distribution of unknown parameter $\theta\in]0,1[$ (i.e. $P_\theta^{X_1}=Bi(1,\theta)$) and, given $X_1=k_1\in\{0,1\}$, $X_2$ has distribution $Bi(1,1-\theta)$ when $k_1=0$ and $Bi(1,\theta)$ when $k_1=1$, i.e. $P_\theta^{X_2|X_1=k_1}=Bi(1,k_1+(1-2k_1)(1-\theta))$. We can think of tossing a coin with probability $\theta$ of getting heads ($=1$) and making a second toss of this coin if it comes up heads on the first toss, or tossing a second coin with probability $1-\theta$ of making heads if the first toss is tails ($=0$). Consider the uniform distribution on $]0,1[$ as the prior distribution $Q$.

So, the joint probability function of $X_1$ and $X_2$ is
\begin{gather*}\begin{split}
f_\theta(k_1,k_2)&=\theta^{k_1}(1-\theta)^{1-k_1}[k_1+(1-2k_1)(1-\theta)]^{k_2}[1-k_1-(1-2k_1)(1-\theta)]^{1-k_2}\\
&=\begin{cases}
\theta(1-\theta)\hbox{ if } k_2=0,\\ (1-\theta)^2\hbox{ if } k_1=0, k_2=1,\\ \theta^2\hbox{ if } k_1=1, k_2=1.
\end{cases}
\end{split}\end{gather*}
It is shown in Nogales(2022b), Example 2, that the Bayes estimator of the conditional probability function $$f_\theta^{X_2|X_1=k_1}(k_2)=[k_1+(1-2k_1)(1-\theta)]^{k_2}[1-k_1-(1-2k_1)(1-\theta)]^{1-k_2}
$$ 
is 
\begin{gather*}
	{f^*_{n,x'}}^{\!\!\!\!X_2|X_1=k_1}(k_2):=\frac{f^*_{n,k'}(k_1,k_2)}{f^*_{n,k',1}(k_1)}=
	\begin{cases}
	\frac{2n+2}{2n+n_{+0}(k')+2n_{01}(k')+3}&\hbox{ if } k_1=k_2=0,\vspace{1ex}\\
	\frac{n_{+0}(k')+2n_{01}(k')+1}{2n+n_{+0}(k')+2n_{01}(k')+3}&\hbox{ if } k_1=0, k_2=1,\vspace{1ex}\\
	\frac{2n+3}{2n+n_{+0}(k')+2n_{01}(k')+4}&\hbox{ if } k_1=1, k_2=0,\vspace{1ex}\\
	\frac{n_{+0}(k')+2n_{01}(k')+1}{2n+n_{+0}(k')+2n_{01}(k')+4}&\hbox{ if } k_1=k_2=1.
		\end{cases}	
\end{gather*}
Theorem 1 shows that this is a strongly consistent estimator of the conditional probability function $f_\theta^{X_2|X_1=k_1}$.
	\end{myexa}

\begin{myexa}\rm 
Let $(X_1,X_2)$ have bivariate normal distribution  
$$N_2\left(
	\left(\begin{array}{c}
		\theta\\
		\theta
		\end{array}\right),\sigma^2
	\left(\begin{array}{cc}
		1 & \rho\\
		\rho & 1
	\end{array}\right)
	\right),
	$$
 and consider the prior distribution $Q=N(\mu,\tau^2)$. 
 It is shown in Nogales (2022b), Example 3, that
 that the conditional distribution 
 $$\left({R^*_{n,x'}}^{\!\!\!\!\!R}\right)^{p_2|p_1=x_1}
 $$
 is $$N\big((1-\rho_1)m_1(x')+\rho_1x_1,\sigma_1^2(1-\rho_1^2)\big)
 $$
 where
\begin{gather*}	
	\rho_1=-\frac{a_n(\rho,\sigma,\tau)+\frac{1-\rho}{1+\rho}}{a_n(\rho,\sigma,\tau)-\frac{1-\rho}{1+\rho}}\cdot\rho,\quad 
	\sigma_1^2=\frac{a_n(\rho,\sigma,\tau)}{a_n(\rho,\sigma,\tau)-\frac{1-\rho}{1+\rho}}\cdot\sigma^2,\\
	m_1(x')=m_2(x')=\frac{s_1(x')+(1+\rho)\frac{\sigma^2}{\tau^2}\mu}{2(1-\rho_1)(1+\rho)^2\sigma^2a_n(\rho,\sigma,\tau)},\\
	a_n(\rho,\sigma,\tau):=2(n+1)(1+\rho)+\frac{\sigma^2}{\tau^2}, \quad s_1(x'):=\sum_i(x'_{i1}+x'_{i2}).
\end{gather*}	
and that its density ${f^*_{n,x'}}^{\!\!\!\!\!p_2|p_1=x_1}$ is the Bayes estimator of the conditional density
$$f_\theta^{X_2|X_1=x_1}(x_2)=\frac1{\sigma\sqrt{2\pi(1-\rho^2)}}
\exp\left\{-\frac{1}{2\sigma^2(1-\rho^2)}[x_2-(1-\rho)\theta-\rho x_1]^2\right\}
$$ 
for the $L^1$-squared loss function. Theorem 1 proves that it is a strongly consistent estimator of the conditional density
$f_\theta^{X_2|X_1=x_1}$.
	\end{myexa}

\section{Appendix.}

Let us briefly recall some basic concepts about Markov kernels, mainly to fix the notations. In the next,  $(\Omega,\mathcal A)$, $(\Omega_1,\mathcal A_1)$ and so on will denote measurable spaces. 

\begin{mydef}\rm  1) (Markov kernel) A Markov kernel
	$M_1:(\Omega,\mathcal A)\pt   (\Omega_1,\mathcal A_1)$ is a map $M_1:\Omega\times\mathcal A_1\rightarrow[0,1]$ such that: 	(i) $\forall \omega\in\Omega$, $M_1(\omega,\cdot)$ is a  probability
	measure on 	$\mathcal A_1$, (ii) $\forall A_1\in\mathcal A_1$, $M_1(\cdot,A_1)$ is $\mathcal A$-measurable.\par
	2)	(Image of a Markov kernel) The image (or {\it probability
		distribution}) of a Markov kernel $M_1:(\Omega,\mathcal A,P)\pt
	(\Omega_1,\mathcal A_1)$ on a probability space is the probability
	measure  $P^{M_1}$ on $\mathcal A_1$ defined by
	$P^{M_1}(A_1):=\int_{\Omega}M_1(\omega,A_1)\,dP(\omega)$.
	\par
	3)  (Composition of Markov kernels) Given two Markov kernels
	$M_1:(\Omega_1,\mathcal A_1)\pt (\Omega_2,\mathcal A_2)$ and
	$M_2:(\Omega_2,\mathcal A_2)\pt (\Omega_3,\mathcal A_3)$, its composition  is defined 	as the Markov kernel $M_2M_1:(\Omega_1,\mathcal A_1)\pt
	(\Omega_3,\mathcal A_3)$ given by
	$$M_2M_1(\omega_1,A_3)=\int_{\Omega_2}M_2(\omega_2,A_3)M_1(\omega_1,d\omega_2).
	$$
\end{mydef}

\begin{myprem}\rm  1) (Markov kernels as extensions of the concept of random variable) The concept of Markov kernel extends the concept of random variable (or measurable map). A random variable $T_1:(\Omega,\mathcal 	A,P)\rightarrow(\Omega_1,\mathcal A_1)$ will be identified with the Markov kernel $M_{T_1}:(\Omega,\mathcal A;P)\pt   (\Omega_1,\mathcal
	A_1)$ defined by $M_{T_1}(\omega,A_1)=\delta_{T_1(\omega)}(A_1)=I_{A_1}(T_1(\omega))$,
	where $\delta_{T_1(\omega)}$ denotes the Dirac measure -the
	degenerate distribution- at the point $T_1(\omega)$, and $I_{A_1}$ is
	the indicator function of the event $A_1$. In particular, the probability distribution $P^{M_{T_1}}$ of $M_{T_1}$ coincides with the probability distribution $P^{T_1}$ of $T_1$ defined as $P^{T_1}(A_1):=P(T_1\in A_1)$\par	
	2) Given a Markov kernel $M_1:(\Omega_1,\mathcal A_1)\pt (\Omega_{2},\mathcal A_{2})$ and a random variable $X_2:(\Omega_2,\mathcal A_2)\rightarrow (\Omega_{3},\mathcal A_{3})$, we have that  $M_{X_2}M_1(\omega_1,A_3)=M_1(\omega_1,X_2^{-1}(A_3))=
	M_1(\omega_1,\cdot)^{X_2}(A_3).$ We write $X_2M_1:=M_{X_2}M_1$.\par
	3) Given a Markov kernel $M_1:(\Omega_1,\mathcal A_1,P_1)\pt (\Omega_{2},\mathcal A_{2})$ we write $P_1\otimes M_1$ for the only probability measure on the product $\sigma$-field $\mathcal A_1\times\mathcal A_2$ such that 
	$$(P_1\otimes M_1)(A_1\times A_2)=\int_{A_1}M_1(\omega_1,A_2)dP_1(\omega_1),\quad A_i\in\mathcal A_i,\, i=1,2.$$

	4) Given two r.v. $X_i:(\Omega,\mathcal A,P)\rightarrow(\Omega_i,\mathcal A_i)$, $i=1,2$, we write $P^{X_2|X_1}$ for the conditional distribution of $X_2$ given $X_1$, i.e. for the Markov kernel $P^{X_2|X_1}:(\Omega_1,\mathcal A_1)\pt(\Omega_2,\mathcal A_2)$ such that $$P^{(X_1,X_2)}(A_1\times A_2)=\int_{A_1}P^{X_2|X_1=x_1}(A_2)dP^{X_1}(x_1),\quad  A_i\in\mathcal A_i,\,i=1,2.	
	$$
	So $P^{(X_1,X_2)}=P^{X_1}\otimes P^{X_2|X_1}$. $\Box$
	
\end{myprem}
	
	Let $(\Omega,\mathcal A,\{P_\theta\colon\theta\in(\Theta,\mathcal T,Q)\})$ be a Bayesian statistical experiment, 
	where $Q$ denotes the prior distribution on the parameter space $(\Theta,\mathcal T)$. We suppose that $P(\theta,A):=P_\theta(A)$ is a Markov kernel $P:(\Theta,\mathcal T)\pt(\Omega,\mathcal A)$. When needed we shall suppose that $P_\theta$ has a density (Radon-Nikodym derivative) $p_\theta$ with respect to a $\sigma$-finite measure $\mu$ on $\mathcal A$ and that the likelihood function $\mathcal L(\omega,\theta):=p_\theta(\omega)$ is $\mathcal A\times\mathcal T$-measurable (this is sufficient to prove that $P$ is a Markov kernel). 
	
	Let $\Pi:=Q\otimes P$, i.e.
	$$\Pi(A\times T)=\int_TP_\theta(A)dQ(\theta), \quad A\in\mathcal A, T\in\mathcal T.
	$$
	The prior predictive distribution is $\beta_Q^*:=\Pi^I$ (the distribution of $I$ with respect to $\Pi$), where $I(\omega,\theta):=\omega$. So
	$$\beta_Q^*(A)=\int_\Theta P_\theta(A)dQ(\theta).
	$$
	The posterior distribution is a Markov kernel $P^*:(\Omega,\mathcal A)\pt(\Theta,\mathcal T)$ such that
	$$\Pi(A\times T)=\int_TP_\theta(A)dQ(\theta)=\int_AP^*_\omega(T)d\beta_Q^*(\omega), \quad A\in\mathcal A, T\in\mathcal T,
	$$
	i.e. such that $\Pi=Q\otimes P=\beta_Q^*\otimes P^*$. This way the Bayesian statistical experiment can be identified with the probability space $(\Omega\times\Theta,\mathcal A\times\mathcal T,\Pi)$, as proposed, for instance, in Florens et al. (1990). 
	
	It is well known that, for $\omega\in\Omega$, the posterior $Q$-density is proportional to the likelihood, i.e.
	$$p^*_{\omega}(\theta):=\frac{dP^*_\omega}{dQ}(\theta)=C(\omega)p_\theta(\omega)
	$$
	where $C(\omega)=[\int_\Theta p_\theta(\omega)dQ(\theta)]^{-1}$.
	
	The posterior predictive distribution on $\mathcal A$ given $\omega$ is 
	$${P_\omega^*}^P(A)=\int_\Theta P_\theta(A)dP_\omega^*(\theta),\quad A\in\mathcal A.
	$$
	This is a Markov kernel 
	$$PP^*(\omega,A):={P_\omega^*}^P(A).$$
	It is readily shown that the posterior predictive density is
	$$\frac{d{P_\omega^*}^P}{d\mu}(\omega')=\int_\Theta p_\theta(\omega')p^*_\omega(\theta)dQ(\theta).
	$$
	
	We know from Nogales (2022a) that
	$$\int_{\Omega\times\Theta}\sup_{A\in\mathcal A}|{P_\omega^*}^P(A)-P_\theta(A)|^2d\Pi(\omega,\theta)\le
	\int_{\Omega\times\Theta}\sup_{A\in\mathcal A}|M(\omega,A)-P_\theta(A)|^2d\Pi(\omega,\theta),$$
	for every Markov kernel $M:(\Omega,\mathcal A)\pt (\Omega,\mathcal A)$ provided that $\mathcal A$ is separable  (recall that a $\sigma$-field is said to be separable, or countably generated, if it contains a countable subfamily which generates it). We also have that, for a real statistic  $X$  with finite mean, the posterior predictive mean
		$$E_{(P_\omega^*)^P}(X)=\int_\Theta\int_\Omega X(\omega')dP_\theta(\omega')dP_\omega^*(\theta)$$
		is the Bayes estimator  of $f(\theta):=E_\theta(X)$, as $E_{(P_\omega^*)^P}(X)=E_{P_\omega^*}(E_\theta(X))$.

\section{Acknowledgements.}
This paper has been supported by the Junta de Extremadura (Spain) under the grant Gr21044.
\vspace{1ex}

\section{References.}

\begin{itemize}

	\item Geisser, S. (1993) Predictive Inference: An Introduction, Springer Science+ Business Media, Dordrecht.
	
	\item Gelman, A., Carlin, J.B., Stern, H.S., Dunson, D.B., Vehtari, A., Rubin, D.B. (2014) Bayesian Data Analysis, 3rd ed., CRC Press.

	\item Ghosal, S., Vaart, A.v.d. (2017) Fundamentals of Noparametric Bayesian Inference, Cambridge University Press, Cambridge UK.

    \item Nogales, A.G. (2022a), On Bayesian Estimation of Densities and Sampling Distributions: the Posterior Predictive Distribution as the Bayes Estimator, Statistica Neerlandica 76(2), 236-250.

    \item Nogales, A.G. (2022b), Optimal Bayesian Estimation of a Regression Curve, a Conditional Density, and a Conditional Distribution, Mathematics 10(8), 1213.

    \item Nogales, A.G. (2022c), On consistency of the Bayes Estimator of the Density, Mathematics, 10(4), 636.

\end{itemize}

\end{document}